\documentclass[11pt, letterpaper]{article}

\newcommand{\bmat}[1]{\begin{bmatrix}#1\end{bmatrix}}
\newcommand{\ie}{\emph{i.e.}}

\title{
Optimal Impact Angle Guidance via First-Order Optimization under Nonconvex Constraints
}

\usepackage{amsmath}
\usepackage{overpic}
\usepackage{rotating}
\usepackage[margin=1in]{geometry}
\newcommand\blfootnote[1]{%
  \begingroup
  \renewcommand\thefootnote{}\footnote{#1}%
  \addtocounter{footnote}{-1}%
  \endgroup
}
\graphicspath{ {./Figures} }

\author{
Gyubin Park\thanks{
Graduate student, Department of Aerospace Engineering, Inha University, South Korea. \texttt{gyubin@inha.edu}
}
\and
Jiwoo Choi\thanks{
Graduate student, Department of Aerospace Engineering, Inha University, South Korea. \texttt{jiwoochoi@inha.edu}
}
\and
Da Hoon Jeong\thanks{
Research engineer, Vehicle Control Development Center, Hyundai Motors Company, South Korea.
}
\and
Jong-Han Kim\thanks{
Associate professor, Department of Aerospace Engineering, Inha University, South Korea. \texttt{jonghank@inha.ac.kr}  (Corresponding author)
}
\blfootnote{
This work was supported in part by Theater Defense Research Center funded by the Korea government(DAPA) under Grant UD200043CD, 
in part by VTOL Technology Research Center for Defense Applications funded by the Korea government(DAPA) under KRIT grant 20-105-E00-005 (2023),
and in part by the KSLV-II Enhancement Program funded by the Korea government(MSIT) (RS-2022-00164702). Contributions are equal.
}
}

\begin{document}

\maketitle

\begin{abstract}

Most of the optimal guidance problems can be formulated as nonconvex optimization problems, which can be solved indirectly by relaxation, convexification, or linearization. Although these methods are guaranteed to converge to the global optimum of the modified problems, the obtained solution may not guarantee global optimality or even the feasibility of the original nonconvex problems. In this paper, we propose a computational optimal guidance approach that directly handles the nonconvex constraints encountered in formulating the guidance problems. The proposed computational guidance approach alternately solves the least squares problems and projects the solution onto nonconvex feasible sets, which rapidly converges to feasible suboptimal solutions or sometimes to the globally optimal solutions. The proposed algorithm is verified via a series of numerical simulations on impact angle guidance problems under state dependent maneuver vector constraints, and it is demonstrated that the proposed algorithm provides superior guidance performance than conventional techniques.

\end{abstract}

\section{Introduction}
\label{sec:introduction}

Most of the optimal guidance and control problems have physical constraints such as the limits on the magnitude or direction of the state variables and control variables, and these constraints are in general expressed in nonconvex form, which in itself has the disadvantage of being difficult to solve efficiently \cite{ref1, ref2}. Classically, several methods have been proposed to deal with the nonconvex constraints indirectly by relaxing or convexifying them \cite{ref3, ref4}, or to solve the optimal control problem for linearized dynamics with regard to the line of sight angle \cite{ref5, ref6, ref13}. 
However, these methods cannot guarantee the optimality of the solution if the linearization assumption around the line of sight vector is invalid due to large divert maneuver. Furthermore, it may be infeasible is the transformed problem does not satisfy the original nonconvex constraints \cite{ref7}. 

In this paper, we propose a computational guidance technique based on a first-order method for handling nonconvex kinematic constraints without any linearization. First-order algorithm converges to the optimal solution by applying separation of variables method to the multiple objective functions and alternately updating each variable, which makes it simple to implement and easy to compute than other second-order algorithms \cite{ref8, ref9}. 

We designed the first-order method that mainly handles the angle constraints between the control vector and the line of sight vector, which can be expressed as dot product of the state variable. The feasible set for this constraint is obviously nonconvex, which is difficult to be handled in general. We show that the orthogonal projection onto this nonconvex set can be solved in explicit form, by transforming the $\ell_2$-norm minimization problem into trigonometric equations. 

The proposed algorithm is designed to iteratively perform the least squares solution and forementioned orthogonal projection onto the nonconvex set. Its validity and effectiveness are verified through a series of numerical simulations on optimal impact angle guidance problems for high-speed target intercept scenarios.

\section{Optimal guidance with state dependent maneuver vector constraints}
\label{sec:optimal_guidance}

\subsection{Problem formulation}

\begin{figure}
\begin{center}
\begin{overpic}[width=0.6\linewidth]{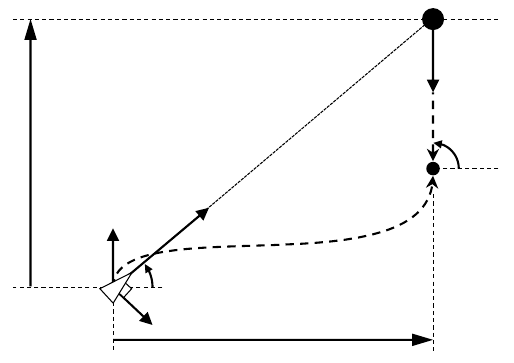}
\put(1,33){\small \begin{turn}{90} Range $(x)$\end{turn}}
\put(44,0){\small Crosstrack $(y)$}
\put(3,9){\small Interceptor}
\put(86,62.5){\small Target}
\put(86,34){\small PIP}
\put(30,7){\small $u_t$}
\put(20.5,28){\small $v_t$}
\put(41,26){\small $l_t$}
\put(89,43){\small $\theta_f$}
\put(31,16){\small $\sigma$}
\put(86,69){\small $\chi_0$}
\put(79.5,51.5){\small $\nu_0$}
\end{overpic}
\end{center}
\caption{Engagement geometry.}
\label{fig1}
\end{figure}
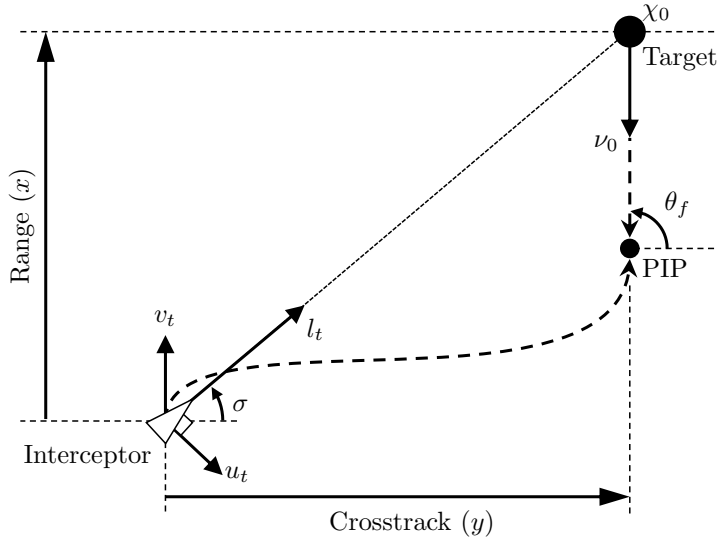

The optimal guidance problem for the engagement geometry shown in Fig. \ref{fig1} is formulated as \eqref{eq1}. 
Assuming an engagement situation where the interceptor largely diverts to the target from its original course, we try the impact angle constrained guidance approach that minimizes the miss distance by forming a head-on trajectory to the target, rather than the commonly used proportional navigation based guidance approach. 
We also assume that a set of lateral thrusters attached perpendicular to the fuselage generate the divert forces for maneuver, while an ideal attitude controller turns its attitude so that its strapdown seeker is always looking towards the target. 
Under these assumptions, the line of sight vector to the target and the body axis always coincide, so that the maneuver acceleration vector is always perpendicular to the line of sight vector.
\begin{subequations}
\label{eq1}
\begin{align}
  \underset{u_0, \dots, u_{N-1}}{\text{minimize}}
                   \quad & \sum_{t=0}^{N-1} \| {u_t} \|^2 \nonumber \\  
                	\text{subject to} \quad & p_{t+1} = p_{t} + \Delta t v_t \label{eq1a}\\
                    \quad & v_{t+1} = v_{t} + \Delta t u_t \label{eq1b}\\
                    \quad & l_t = \chi_0 + \Delta t \nu_0 t - p_t \label{eq1c}\\
                    \quad & l_N = 0 \label{eq1d}\\
                    \quad & \angle v_N = \theta_f \label{eq1e}\\
                    \quad & \angle (u_t, l_t) = \pi/2 \label{eq1f}\\
                    \quad & \| u_t \| \leq u_\text{ub} \label{eq1g}
\end{align}
\end{subequations}

The horizon size and the step size is denoted by $N$ and $\Delta t$, and the maneuver acceleration, position, and velocity of the missile at time $t$ are denoted by $u_t$, $p_t$, and $v_t$, respectively.
Equations \eqref{eq1a} and \eqref{eq1b} represent the dynamics of the missile, and \eqref{eq1c} describes the line of sight vector to the constant velocity target with velocity $\nu_0$ from the initial position $\chi_0$. 
Equation \eqref{eq1d} and \eqref{eq1e} are the terminal position and the terminal impact angle constraints. 
Equation \eqref{eq1f} and \eqref{eq1g} states that the limited maneuver acceleration is always perpendicular to the line of sight vector. The dynamic constraints \eqref{eq1a} and \eqref{eq1b} and the terminal constraints \eqref{eq1d} and \eqref{eq1e} can be encoded into a single linear equation $A u - b = \bmat{l_N^T,~\tan \theta_f - v_{y, N} / v_{x, N}}^T = 0$, where we used the stack notation $u=(u_0,~\cdots,~u_{N-1})=\bmat{u_0^T & \cdots & u_{N-1}^T}^T$. This results in \eqref{eq2} where $v_{x,t}$ and $v_{y,t}$ stand for the $x$- and $y$- components of $v_t$.
\begin{subequations}
\label{eq2}
\begin{align}
  \underset{u}{\text{minimize}}
                    \quad & \| u \|^2 \nonumber \\
  \text{subject to} \quad & A u - b = 0 \label{eq2a} \\
                    \quad & \angle (u_t, l_t) = \pi/2 \label{eq2b} \\
                    \quad & v_{y, N} \ge 0 \label{eq2c} \\
                    \quad & \| u_t \| \leq u_\text{ub} \label{eq2d}
\end{align}
\end{subequations}

Note that the terminal impact angle constraint was implemented by using the tangent function \eqref{eq3a} as follows with a halfspace constraint \eqref{eq3b}, in order for avoiding the ambiguity originated from the periodicity of the tangent function. Also note that the equation in \eqref{eq3a} can be encoded into \eqref{eq2a}.
\begin{subequations}
\label{eq3}
\begin{align}
\frac{v_{y, N}}{v_{x, N}} &= \tan{\theta_f}\label{eq3a} \\
v_{y, N} &\ge 0  \label{eq3b}
\end{align}
\end{subequations}

Since the orthogonality condition \eqref{eq2b} is nonconvex, we will use a first-order method that handles the constraints by computing the orthogonal projection onto the feasible set. In this paper, we propose an efficient algorithm for computing the orthogonal projection onto the nonconvex set described by the orthogonality condition between state vectors, given in \eqref{eq2b}.

\subsection{First-order methods}
First-order methods are optimization algorithms that use the gradient or subgradient of the objective function.
The methods are more robust, easier to implement, and more efficient in memory usage than the second order methods using the Hessian information. The second order methods typically converge in fewer iterations, however requires heavier computational complexity for each iteration~\cite{ref8}.

In this paper, we perform optimization using a first-order method, the alternating direction method of multipliers (ADMM). The ADMM is a primal-dual algorithm combining dual ascent with the method of multipliers. It sequentially updates the primal variables, followed by dual ascent \cite{ref8}.

A general optimization problem with an objective function $f(\cdot)$ and a feasible set $\mathcal{C}$ can be expressed as
\begin{equation}
\label{eq4}
    \begin{aligned}
    \text{minimize}\quad   & f(x)\\
    \text{subject to}\quad & x \in \mathcal{C}
\end{aligned}
\end{equation}
which can be reformulated into an ADMM form as follows.
\begin{equation}
\label{eq5}
    \begin{aligned}
    \text{minimize}\quad & f(x) + I_{\mathcal{C}}(z)\\
    \text{subject to}\quad & x =z
\end{aligned}
\end{equation}

Here, $I_\mathcal{C}(\cdot)$ is defined as the indicator function that returns $0$ or $\infty$ depending on whether the input is in the feasible set or not.


In general when the function $f$ and the set $\mathcal{C}$ are convex, we are guaranteed to find the globally optimal solution for the problem in \eqref{eq4}. This is due to nonexpansivity of the proximal operator of convex functions. If either the function $f$ or the set $\mathcal{C}$ is nonconvex, the global convergence is not guaranteed, however we can frequently find practically good solutions in many cases~\cite{ref10}. In addition, global convergence can be guaranteed in some nonconvex cases~\cite{ref11}. In this study, we consider the case where the set $\mathcal{C}$ is nonconvex.

\subsection{First-order method for optimal impact angle guidance}

The primal variables for the nonconvex constraint and the inequality constraints can be written as,
\begin{subequations}
\label{eq7}
\begin{align}
  z_1 &=
  u
  \in
  \mathcal{C}_1 \label{eq7a} \\
  z_2 &=
  v_{y, N} =
  P u - q
  \in
  \mathcal{C}_2 \label{eq7b} \\
  z_3 &=
  \begin{bmatrix}
    u^T & l^T
  \end{bmatrix}^T
  = Gu - h
  \in
  \mathcal{C}_3 \label{eq7c}
\end{align}
\end{subequations}
with the feasible set defined by,
\begin{subequations}
\label{eq8}
\begin{align}
\mathcal{C}_1(u) &= \{ u \ \vert\  \|u_t\|\le u_\text{ub} ,\ \forall t \} \label{eq8a} \\
\mathcal{C}_2(v) &= \{ v \ \vert\  v_{y,N} \ge 0 ,\ \forall t \} \label{eq8b} \\
\mathcal{C}_3(u,l) &= \{ (u,l) \ \vert\  \angle (u_t, l_t) = \pi/2 ,\ \forall t \} \label{eq8c}
\end{align}
\end{subequations}
where we used a stack notation $(u,~l) = \bmat{u^T & l^T}^T$.
For nonconvex constraint $\angle (u_t, l_t) = \pi/2$, we first need the position of the missile throughout the entire time horizon, which leads to time series of the line of sight vector.
Here, $z_3$ is the primal variable consisting of the control input and the line of sight vector. These are projected onto their own feasible set, which in this case is $\mathcal{C}_3$. 

Other primal variables $z_1$ and $z_2$ represent inequality constraints, which are projected onto feasible set $\mathcal{C}_1$ and $\mathcal{C}_2$, respectively.
Taken together, the ADMM form of the original problem can be written as follows.
\begin{equation}
\label{eq9}
\begin{aligned}
  \underset{u, z_1, z_2, z_3}{\text{minimize}}
                    \quad & \| u \|^2 + I_{\mathcal{C}_1}(z_1) + I_{\mathcal{C}_2}(z_2) + I_{\mathcal{C}_3}(z_3) \\
  \text{subject to} \quad & A u - b = 0 \\
                    \quad & u = z_1 \in \mathcal{C}_1 \\
                    \quad & P u - q = z_2 \in \mathcal{C}_2 \\
                    \quad & G u - h = z_3 \in \mathcal{C}_3
\end{aligned}
\end{equation}

With the intersection of all the feasible set defined by $\mathcal{C} = \mathcal{C}_1 \cap \mathcal{C}_2 \cap \mathcal{C}_3$, we can write the augmented Lagrangian as,
\begin{equation}
\label{eq10}
\begin{aligned}
L_\rho (u, z, w) &= 
\| u \|^2
+ I_\mathcal{C}(z) + \frac{\rho}{2}
\|
    M u - n - z + w
\|^2 \\
\text{where }
M
=&
\begin{bmatrix}
  A \\
  I \\
  P \\
  G
\end{bmatrix},~
n
=
\begin{bmatrix}
  b \\
  0 \\
  q \\
  h
\end{bmatrix},~
z
=
\begin{bmatrix}
  0 \\
  z_1 \\
  z_2 \\
  z_3
\end{bmatrix},~
w
=
\begin{bmatrix}
  w_0 \\
  w_1 \\
  w_2 \\
  w_3
\end{bmatrix}
\end{aligned}
\end{equation}
and the optimal solution can be found by iteratively computing \eqref{eq11}-\eqref{eq15}, which are the update processes for the set of primal variables $u$, $z_1$, $z_2$, $z_3$ and the scaled dual variable $w$. 

\begin{itemize}

\item $u$-update step (Least squares solution):
\begin{equation}
\label{eq11}
\begin{aligned}
    u^{(k+1)}
    &=
    \underset{u}{\text{argmin}}
    \left\{
         \| u \|^2 + \displaystyle\frac{\rho}{2}
         \| M u - n - z^{(k)} + w^{(k)} \|^2
    \right\}\\
    &=
    \frac{\rho}{2}
    \left(
      I + \displaystyle\frac{\rho}{2}M^TM
    \right)^{-1}
    M^T
    \left(
      n + z^{(k)} - w^{(k)}
    \right)
\end{aligned}
\end{equation}

\item $z_1$-update step (Projection onto $\mathcal{C}_1$):
\begin{equation}
\label{eq12}
\begin{aligned}
    z_1^{(k+1)}
    &=
    \underset{z_1}{\text{argmin}}
    \left\{
         I_{\mathcal{C}_1}(z_1) + \displaystyle\frac{\rho}{2}
         \| u^{(k+1)} - z_1 + w_1^{(k)} \|^2
    \right\} \\
    &=
    \Pi_{\mathcal{C}_1}
    \left(
      u^{(k+1)} + w_1^{(k)}
    \right)
\end{aligned}
\end{equation}

\item $z_2$-update step (Projection onto $\mathcal{C}_2$):
\begin{equation}
\label{eq13}
\begin{aligned}
    z_2^{(k+1)}
    &=
    \underset{z_2}{\text{argmin}}
    \bigg{\{}
         I_{\mathcal{C}_2}(z_2) \\
         &\qquad\qquad\quad + \displaystyle\frac{\rho}{2}
         \| Pu^{(k+1)} - q - z_2 + w_2^{(k)} \|^2
    \bigg{\}} \\
    &=
    \Pi_{\mathcal{C}_2}
    \left(
      Pu^{(k+1)} - q + w_2^{(k)}
    \right)
\end{aligned}
\end{equation}

\item $z_3$-update step (Projection onto $\mathcal{C}_3$):
\begin{equation}
\label{eq14}
\begin{aligned}
    z_3^{(k+1)}
    &=
    \underset{z_3}{\text{argmin}}
    \bigg{\{}
         I_{\mathcal{C}_3}(z_3) \\
         &\qquad\qquad\quad + \displaystyle\frac{\rho}{2}
         \| Gu^{(k+1)} - h - z_3 + w_3^{(k)} \|^2
    \bigg{\}} \\
    &=
    \Pi_{\mathcal{C}_3}
    \left(
      Gu^{(k+1)} - h + w_3^{(k)}
    \right)
\end{aligned}
\end{equation}

\item $w$-update step (Dual ascent):
\begin{equation}
\label{eq15}
w^{(k+1)} = w^{(k)} + Mu^{(k+1)} - n - z^{(k+1)}
\end{equation}

\end{itemize}

The $u$-update step in \eqref{eq11} is a multi-objective weighted least squares problem which is easy to solve. 
Computing the inverse matrix, which can be computationally expensive, does not change with iterations, therefore it needs to be calculated only once when the algorithm is initialized. 
Moreover, most of the appearing matrices, $M$, $n$, and so on, are sparse, which helps us to greatly reduce the computational complexity.

Due to the inherent sensitivity of procedures~\eqref{eq11} to \eqref{eq15} when applied to nonconvex problem instances, particularly to initial conditions, hyperparameters, and the sequence of variable updates as highlighted in prior literature~\cite{ref11, ref12}, a strategic ordering of updates is essential for improved convergence. Consequently, the $z_1$-update and $z_2$-update steps, involving projections onto convex sets, are prioritized over the $z_3$-update step, which entails a projection onto a nonconvex set that can be expansive. This sequencing of variable updates aims to enhance convergence characteristics.

The $z_1$-update step in \eqref{eq12} shrinks $\|u_t\|$ to $u_\text{ub}$ when $\|u_t\| > u_\text{ub}$, and the $z_2$-update step in \eqref{eq13} returns the positive part only, \ie, $(v_{y, N})_+$.
\begin{equation}
\label{eq16}
\begin{aligned}
    \Pi_{\mathcal{C}_1}
    &\left(
      u_t^{(k+1)} + w_1^{(k)}
    \right) \\
    &=
    \begin{cases}
      u_t^{(k+1)} + w_{1, t}^{(k)}, & \text{if } \| u_t^{(k+1)} + w_{1, t}^{(k)} \| \leq u_\text{ub} \\
      u_\text{ub}
      \displaystyle
      \frac{u_t^{(k+1)} + w_{1, t}^{(k)}}{\|u_t^{(k+1)} + w_{1, t}^{(k)}\|}, & \text{otherwise}
    \end{cases}
\end{aligned}
\end{equation}
\begin{equation}
\label{eq17}
\Pi_{\mathcal{C}_2}
\left(
  Pu^{(k+1)} - q + w_2^{(k)}
\right)
= 
\left( Pu^{(k+1)} - q + w_2^{(k)} \right)_+
\end{equation}

Solving the Euclidean projection from the $z_3$-update step of \eqref{eq14} to the nonconvex feasible set is the core of the algorithm proposed in this paper. We will show in Section~\ref{sec:projection} that the projection result can be obtained in closed form.

The dual update step in \eqref{eq15} is also simple matrix-vector computations, which is similarly easy to compute.

\begin{figure}
\begin{center}
    \begin{overpic}[width=\linewidth]{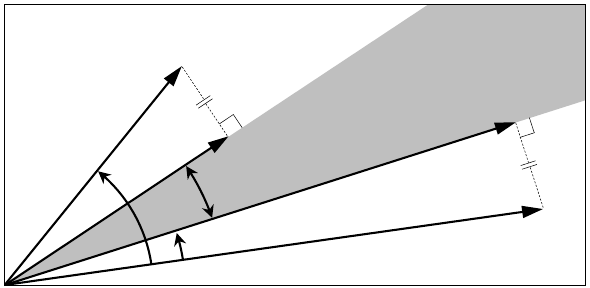}
        \put(65,32){\small $\mathcal{C}_\text{narr}$}
        \put(89,9){\small $\alpha_t$}
        \put(26,38){\small $\beta_t$}
        \put(83,30){\small $z_{\alpha, t}$}
        \put(38,22){\small $z_{\beta, t}$}
        \put(24,11){\small $\phi$}
        \put(35,16){\small $\theta$}
        \put(32,6.5){\small $\psi$}
    \end{overpic}
\end{center}
\vspace{-0.2cm}
\caption{Projection onto the interior of a state dependent cone with apex angle $\theta$, defined by
 $\mathcal{C}_\text{narr} = \{ (\alpha,\beta)\ \vert\ \angle(\alpha,\beta)\le\theta \}$.
}
\label{fig2}
\end{figure}

\begin{figure}
\begin{center}
    \begin{overpic}[width=\linewidth]{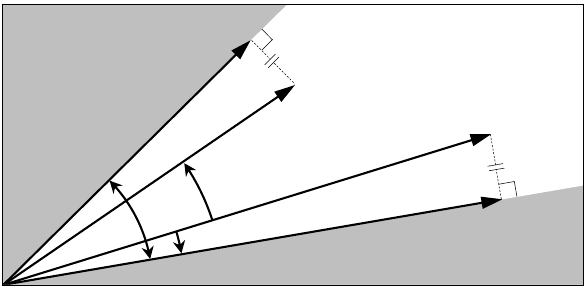}
        \put(13,32){\small $\mathcal{C}_\text{wide}$}
        \put(70,5){\small $\mathcal{C}_\text{wide}$}
        \put(80,28){\small $\alpha_t$}
        \put(48,29){\small $\beta_t$}
        \put(82,11){\small $z_{\alpha, t}$}
        \put(35,43){\small $z_{\beta, t}$}
        \put(35,16){\small $\phi$}
        \put(24,11){\small $\theta$}
        \put(31.5,7.5){\small $\psi$}
    \end{overpic}
\end{center}
\vspace{-0.2cm}
\caption{Projection onto the exterior of a state dependent cone with apex angle $\theta$, defined by
 $\mathcal{C}_\text{wide} = \{ (\alpha,\beta)\ \vert\ \angle(\alpha,\beta)\ge\theta \}$.
}
\label{fig3}
\end{figure}

\section{Projection onto the nonconvex feasible set for angular constraints}
\label{sec:projection}

In this section, we present a closed form solution for finding the orthogonal projection onto the set of vectors whose inbetween angle is equal to $\theta$; this corresponds to finding the orthogonal projection onto the surface of a state dependent cone with apex angle $\theta$. 
Specifically, note that we require $\theta=\pi/2$ in this paper.
\begin{equation}
\label{eq18}
	\mathcal{C}_3(\alpha,\beta) = \left\{ (\alpha,\beta) \ \vert\ \angle(\alpha, \beta) = \theta \right\}
\end{equation}

Note that the symmetric axis of the cone onto which the variables are projected is not fixed in space, instead it should be simultaneously found with the orthogonal projection. 
We also note that the feasible set $\mathcal{C}_3$ can be seen as the intersection of the following two sets, therefore the projection can be separately computed. In other words, the points can be projected onto $\mathcal{C}_\text{narr}$ when the angle between them is larger than $\theta$, or they can be projected onto $\mathcal{C}_\text{wide}$ otherwise. See Figure~\ref{fig2} and Figure~\ref{fig3}.
\begin{equation}
\label{eq19}
\begin{aligned}
	\mathcal{C}_\text{narr}(\alpha,\beta) &= \left\{ (\alpha,\beta) \ \vert\ \angle(\alpha, \beta) \le \theta \right\} \\
	\mathcal{C}_\text{wide}(\alpha,\beta) &= \left\{ (\alpha,\beta) \ \vert\ \angle(\alpha, \beta) \ge 	\theta \right\}
\end{aligned}
\end{equation}

More specifically, we can express the $z$-update step as follows with 
$\eta^{(k+1)} = Gu^{(k+1)} - h + w_3^{(k)}$,
\begin{equation}
\label{eq20}
\begin{aligned}
    z_3^{(k+1)}
    &=
    \underset{z_3}{\text{argmin}}
    \left\{
         I_{\mathcal{C}_3}(z_3) + \displaystyle\frac{\rho}{2}
         \| Gu^{(k+1)} - h - z_3 + w_3^{(k)} \|^2
    \right\} \\
    &=
    \underset{z_3}{\text{argmin}}
    \left\{
         I_{\mathcal{C}_3}(z_3) + \displaystyle\frac{\rho}{2}
         \| \eta^{(k+1)} - z_3 \|^2
    \right\} \\
    &=
    \begin{cases}
      \eta^{(k+1)}, & \text{if }\eta^{(k+1)} \in \mathcal{C}_3 \\
      \underset{z_3}{\text{argmin}}\
      \|
           \eta^{(k+1)} - z_3 
      \|^2, & \text{otherwise}
    \end{cases}
\end{aligned}
\end{equation}
and 
\begin{equation}
\label{eq21}
\begin{aligned}
\eta^{(k+1)} &=
\begin{bmatrix}
  u^{(k+1)} + w_\alpha^{(k)} \\
  l^{(k+1)} + w_\beta^{(k)}
\end{bmatrix}
=
\begin{bmatrix}
  \alpha \\ \beta
\end{bmatrix} \\ 
z_3^{(k+1)} &=
\begin{bmatrix}
  z_{\alpha} \\
  z_{\beta}
\end{bmatrix}. 
\end{aligned}
\end{equation}

Given $\angle\begin{pmatrix}\alpha_t, \beta_t\end{pmatrix}=\phi$ and $\angle\begin{pmatrix}z_{\alpha, t}, z_{\beta, t}\end{pmatrix}=\theta$, we can find the projection onto the sets $\mathcal{C}_\text{narr}$ and $\mathcal{C}_\text{wide}$ by finding the angle $\psi$ that minimizes the sum of the squared distances, $\|z_{\alpha, t} - \alpha_t\|^2+ \|z_{\beta, t} - \beta_t\|^2$.
\begin{subequations}
\label{eq22}
\begin{align}
    z_{\alpha, t}, z_{\beta, t}
    = &
    \underset{z_{\alpha, t}, z_{\beta, t}}{\text{argmin}}
    \left\{
         \| \alpha_t - z_{\alpha, t} \|^2
         +
         \| \beta_t - z_{\beta, t} \|^2
    \right\} \label{eq22a} \\
    & \qquad\qquad\quad~~ \Updownarrow \nonumber \\
    \psi
    =
    \underset{\psi}{\text{argmin}} &
    \left\{
         \| \alpha_t \|^2
         \sin^2{\psi}
         +
         \| \beta_t \|^2
         \sin^2{(\phi-\theta-\psi)}
    \right\} \nonumber \\
    =
    \underset{\psi}{\text{argmin}} &
    \ g (\psi) \label{eq22b}
\end{align}
\end{subequations}

If $\pi > \phi > \theta > 0$, the projection onto $\mathcal{C}_\text{narr}$ can be written as \eqref{eq22}, and the function $g(\psi)$ is unimodal with the minimizer $\psi \in \bmat{0,~\phi - \theta}$ only if $\phi - \theta \leq \pi/2$.
\begin{equation}
\label{eq23}
\begin{aligned}
  \gamma
  & =
  \| \alpha_t \|^2
  +
  \| \beta_t \|^2
  \cos
  \left(
    2 ( \phi - \theta )
  \right) \\
  \delta
  & =
  \| \beta_t \|^2
  \sin 
  \left(
    2 ( \phi - \theta )
  \right) \\
  \Omega
  & =
  \tan^{-1}
  \left(
    \displaystyle
    {\delta}/{\gamma}
  \right) \\
\end{aligned}
\end{equation}

Since $g(\psi)$ is now unimodal, we can isolate $\psi$ through some trigonometric identities, which satisfies $\nabla g = 0$.
\begin{equation}
\label{eq24}
\begin{aligned}
  \nabla g
  & =
  \gamma  \sin ( 2 \psi )
  -
  \delta \cos ( 2 \psi )\\
  & =
  \sqrt{
    {\gamma}^2 + {\delta}^2
  }
  \sin ( 2\psi - \Omega ) \\
  & =
  0
\end{aligned}
\end{equation}

As a result, the minimizer $\psi$ can be uniquely determined as follows.
\begin{equation}
\label{eq25}
\psi = {\Omega}/{2}
\end{equation}

We can write the dot product of $\alpha_t$, $\beta_t$ and $z_{\alpha, t}$, $z_{\beta, t}$ in matrix form, to find the projection $z_{\alpha, t}$ and $z_{\beta, t}$ by solving the linear equation below,
\begin{equation}
\label{eq26}
\begin{bmatrix}
  \alpha_t^T & 0 \\
  \beta_t^T & 0 \\
  0 & \alpha_t^T \\
  0 & \beta_t^T \\
\end{bmatrix}
\begin{bmatrix}
  z_{\alpha, t} \\
  z_{\beta, t}
\end{bmatrix}
=
d_\text{narr, wide}
\end{equation}
where the right hand side, $d_\text{narr, wide}$, can be either $d_\text{narr}$ or $d_\text{wide}$ depending on which set among $\mathcal{C}_\text{narr}$ and $\mathcal{C}_\text{wide}$ we compute the projection onto.
\begin{itemize}
\item Projection onto $\mathcal{C}_\text{narr}$:
\begin{equation}
\label{eq27}
d_\text{narr}
=
\begin{bmatrix}
  \|\alpha_t\|\|\alpha_t\| \cos^2\psi \\
  \|\alpha_t\|\|\beta_t\| \cos\psi \cos(\phi-\psi) \\
  \|\alpha_t\|\|\beta_t\| \cos(\phi-\theta-\psi) \cos(\theta+\psi) \\
  \|\beta_t\|\|\beta_t\| \cos^2(\phi-\theta-\psi)
\end{bmatrix}
\end{equation}

\item Projection onto $\mathcal{C}_\text{wide}$:
\begin{equation}
\label{eq28}
d_\text{wide}
=
\begin{bmatrix}
  \|\alpha_t\|\|\alpha_t\| \cos^2\psi \\
  \|\alpha_t\|\|\beta_t\| \cos\psi \cos(\phi+\psi) \\
  \|\alpha_t\|\|\beta_t\| \cos(\theta-\phi-\psi) \cos(\theta-\psi) \\
  \|\beta_t\|\|\beta_t\| \cos^2(\theta-\phi-\psi)
\end{bmatrix}
\end{equation}
\end{itemize}

Note that we can easily extend the results to three dimensional space, where we can find $z_{\alpha, t}$ and $z_{\beta, t}$ on the plane containing $\alpha_t$, $\beta_t$ described by the normal vector $\alpha_t \times \beta_t$. In Section~\ref{sec:numerical_examples}, we will validate the proposed algorithm via a series of numerical examples and compare the guidance performance with the classical impact angle guidance techniques.

\begin{figure}
\begin{center}
    \begin{overpic}[width=\linewidth]{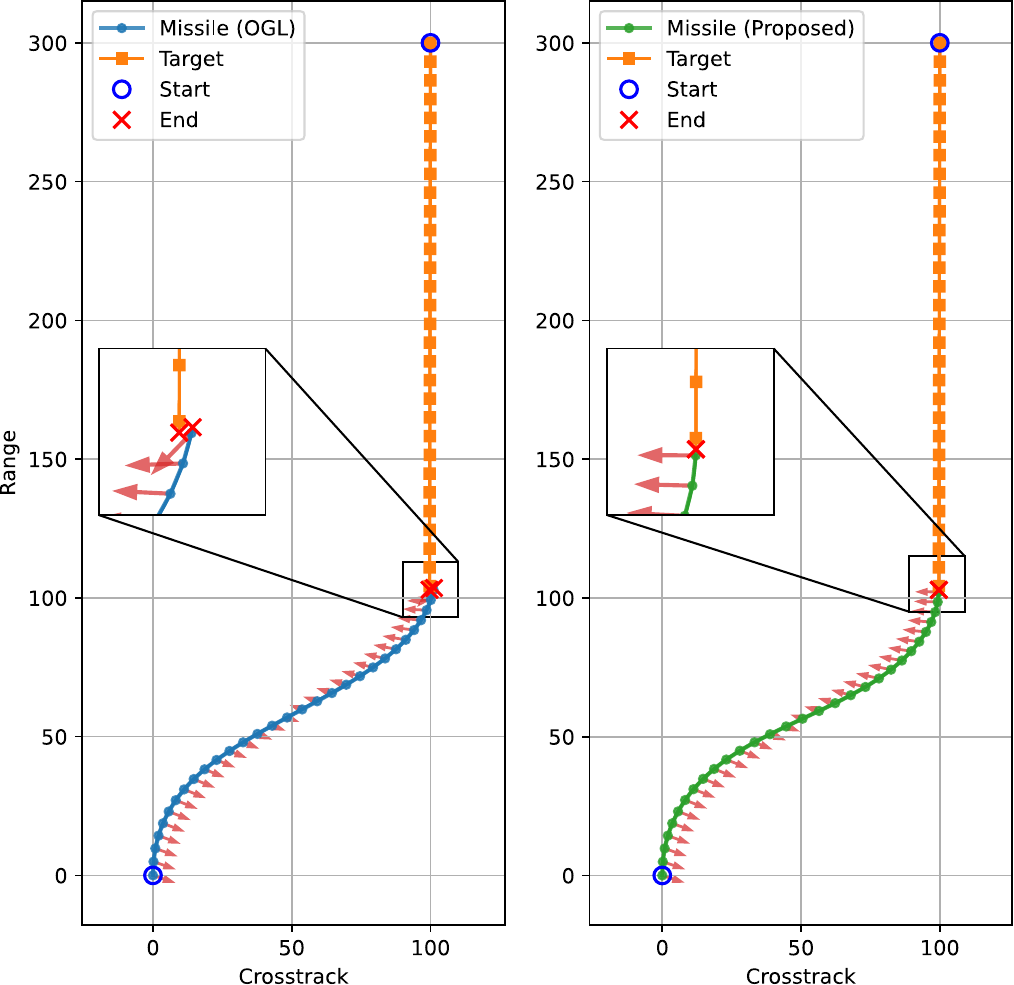}
    \end{overpic}
\end{center}
\caption{Engagement trajectories and maneuver acceleration commands.}
\label{fig4}
\end{figure}

\section{Numerical examples}
\label{sec:numerical_examples}

\subsection{Classical impact angle guidance}

In order to compare the performance of the optimal computational guidance solution obtained in Sections~\ref{sec:optimal_guidance} and \ref{sec:projection}, we consider the classical impact angle guidance solution as a reference. The classical impact angle guidance technique, often referred to as the optimal guidance law (OGL), can be derived from a continuous-time LQR problem linearized around the initial line of sight vector as follows \cite{ref13},
\begin{equation}
\label{eq29}
a_\text{OGL} = V_c\left(4\dot{\sigma}+2\displaystyle\frac{\sigma-\theta_f}{t_\text{go}}\right)
\end{equation}
where $V_c$ is the closing velocity between the interceptor and the target. The line of sight angle and the line of sight rate are denoted by, $\sigma$ and $\dot{\sigma}$ , respectively. The desired terminal impact angle and the time-to-go until impact are written by $\theta_f$ and $t_\text{go}$. 

This solution is very robust and easy to implement, however the optimality easily breaks down in the large divert scenarios where the vehicle moves significantly from the initial course to which the problem is linearized.

\subsection{Simulation results}

\begin{figure}
\begin{center}
    \begin{overpic}[width=0.95\linewidth]{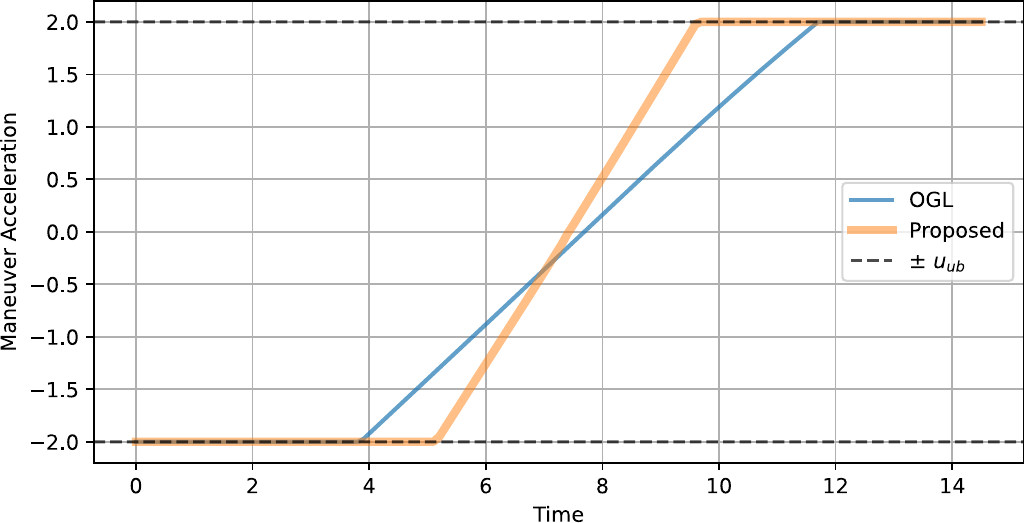}
    \end{overpic}
\end{center}
\vspace{-0.2cm}
\caption{Maneuver acceleration and acceleration limit.}
\label{fig5}
\end{figure}

\begin{figure}
\vspace{0.5cm}
\begin{center}
    \begin{overpic}[width=0.95\linewidth]{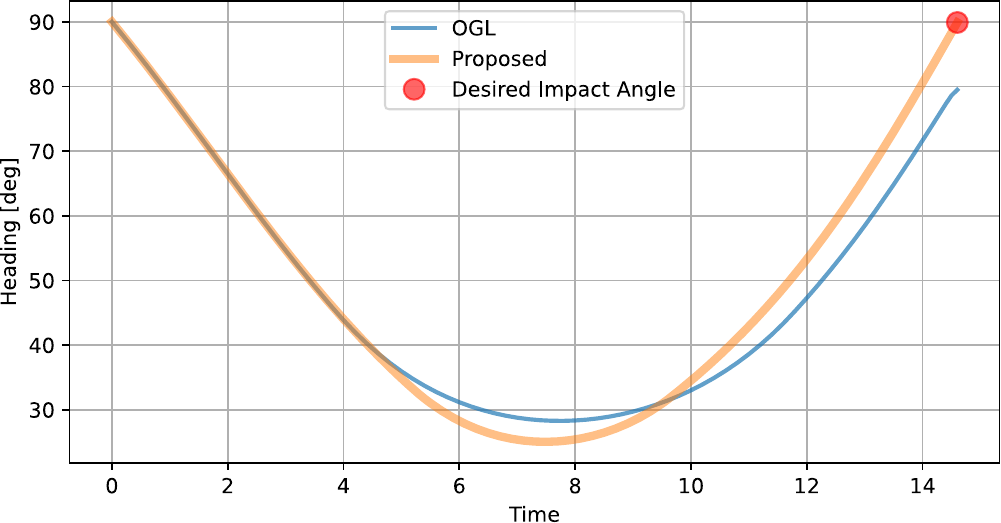}
    \end{overpic}
\end{center}
\vspace{-0.2cm}
\caption{Heading angle and desired impact angle.}
\label{fig6}
\end{figure}

\begin{figure}
\begin{center}
    \begin{overpic}[width=0.95\linewidth]{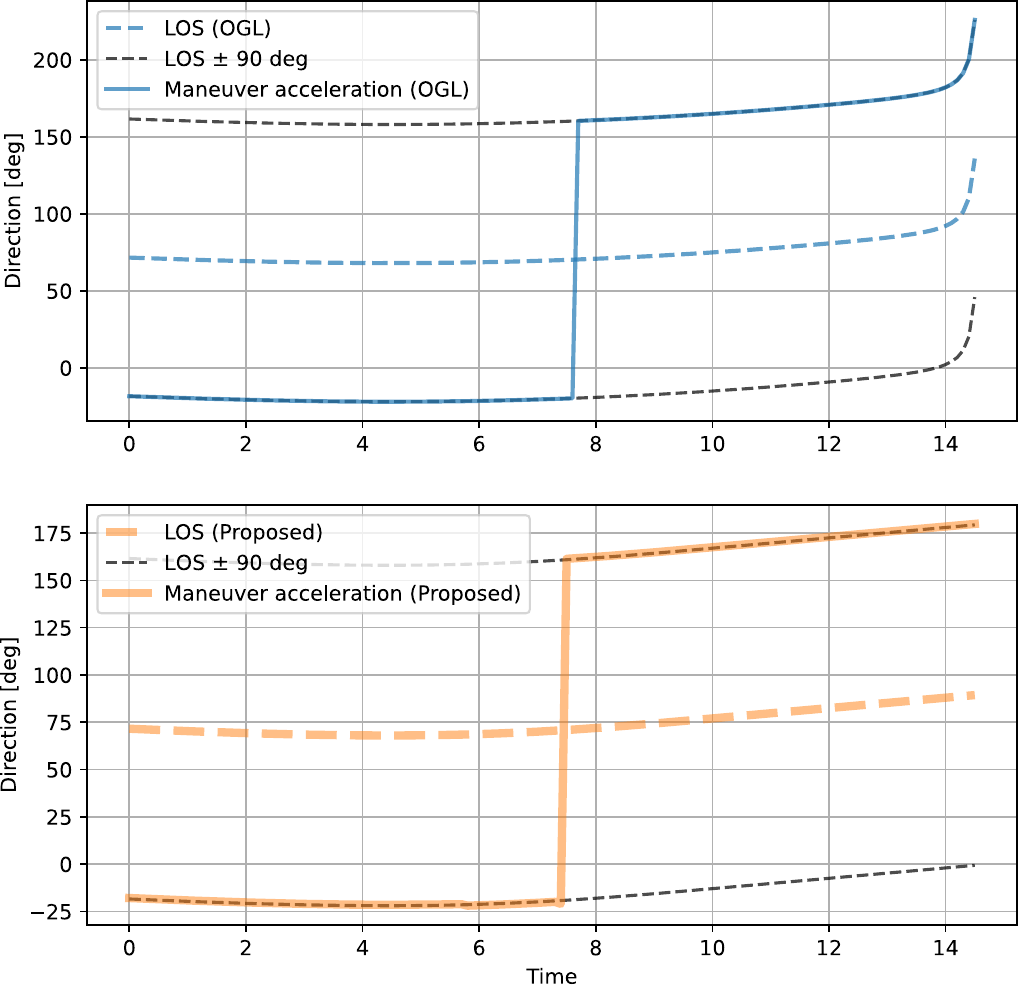}
    \end{overpic}
\end{center}
\vspace{-0.2cm}
\caption{Direction of the line of sight (LOS) vector and the maneuver acceleration vector.}
\label{fig7}
\end{figure}

\begin{figure}
\vspace{0.5cm}
\begin{center}
    \begin{overpic}[width=0.95\linewidth]{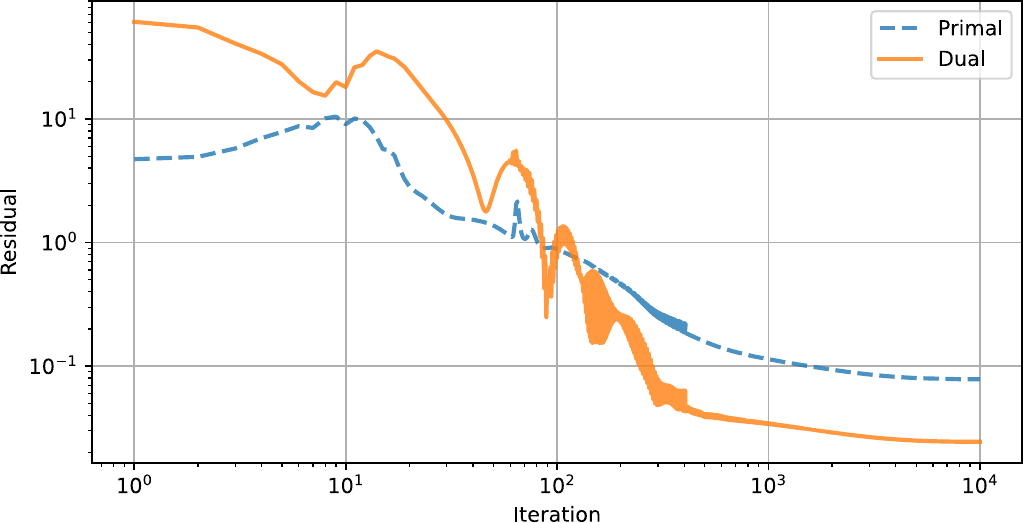}
    \end{overpic}
\end{center}
\vspace{-0.2cm}
\caption{Primal residual and dual residual.}
\label{fig8}
\end{figure}

\begin{figure}
\begin{center}
    \begin{overpic}[width=\linewidth]{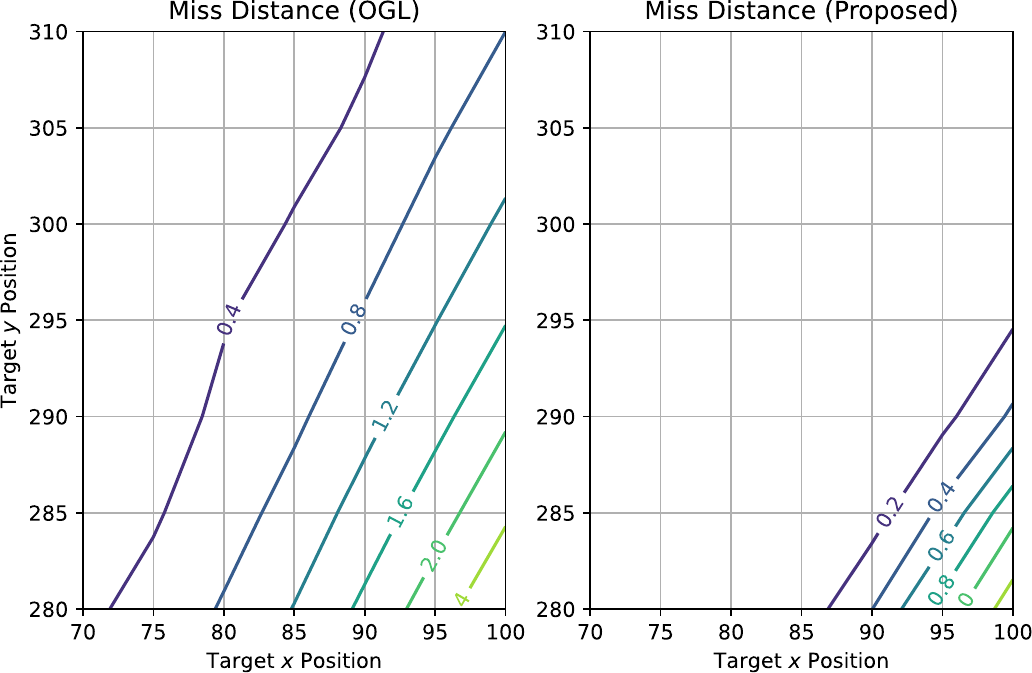}
    \end{overpic}
\end{center}
\vspace{-0.2cm}
\caption{Miss distance for various initial positions of the target.}
\label{fig9}
\end{figure}

\begin{figure}[t]
\vspace{0.5cm}
\begin{center}
    \begin{overpic}[width=\linewidth]{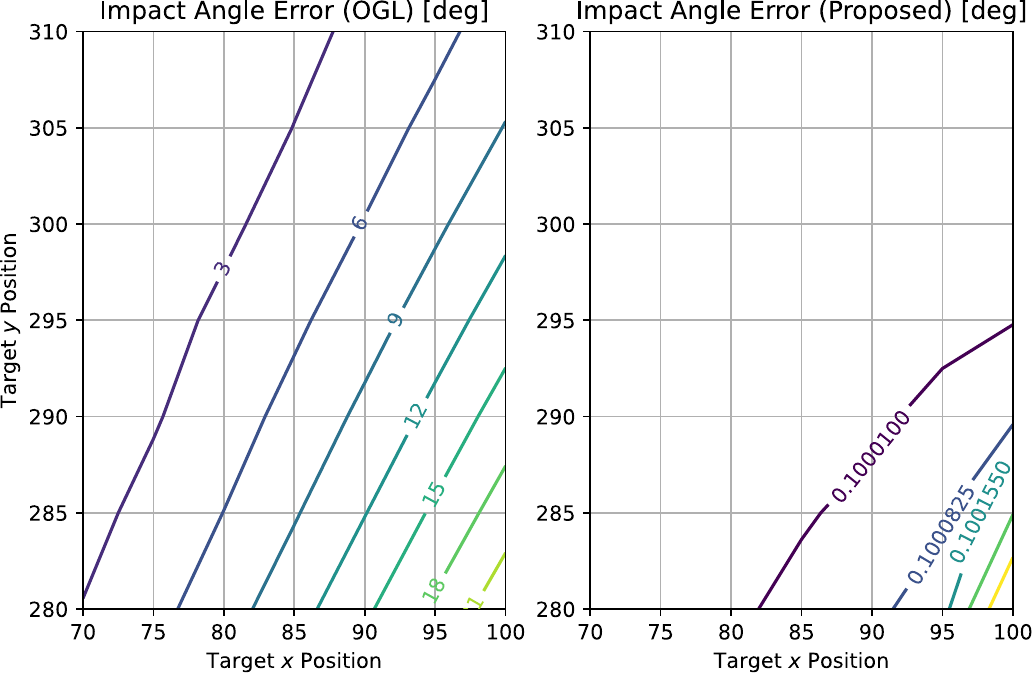}
    \end{overpic}
\end{center}
\vspace{-0.2cm}
\caption{Impact angle error for various initial positions of the target.}
\label{fig10}
\end{figure}

The proposed technique and the classical OGL were configured for an engagement scenario with the terminal impact angle constraint of $\theta_f=\pi/2$ and the initial positions shown at Figure \ref{fig4}, and the results are shown from Figure \ref{fig4} to Figure \ref{fig10}.

Here, the horizon size $N$ required for the proposed method is set equal to the horizon size in the OGL engagement scenario to facilitate fair comparison between the two algorithms. By the way, the minimum $N$ that yields a feasible solution can be found by bisection over $N$, if solving for the minimum time control under given constraints is required.

Figure \ref{fig4} shows the engagement trajectories for both algorithms, and the arrows shown are the maneuver acceleration as control input at each time step. We can observe that the proposed algorithm successfully forms a head-on trajectory with the target and intercepts the target, while the OGL fails to form a head-on trajectory and further fails to intercept the target. This can also be seen from Figure \ref{fig6}.

Figure \ref{fig5} and Figure \ref{fig6} show the magnitude of the maneuver acceleration and the heading angle for each algorithm. 
In the case of the OGL solution, due to the acceleration limit at the beginning of the engagement (when the $\sigma-\theta_f$ is large) and just before the impact (when $\dot{\sigma}$ diverges), the desired terminal impact angle is not achieved and it resulted in large miss distance. On the other hand, it can be seen that the proposed approach uses $u_t$ significantly throughout the entire engagement to achieve the desired terminal impact angle.

Figure \ref{fig7} shows the direction of the line of sight vector and the direction of the maneuver acceleration vector, which validates the projection operation introduced in Section~\ref{sec:projection}. For the OGL results on the top, the angle difference is maintained at $\pi/2$ because $a_\text{OGL}$ calculated from equation \eqref{eq29} is applied to the direction perpendicular to the line of sight vector, and the proposed algorithm on the bottom also confirms that the projection operation is correctly performed and satisfies the orthogonality constraint.

Figure \ref{fig8} shows the primal residual $r^{(k+1)}$ and the dual residual $s^{(k+1)}$ calculated from iterating \eqref{eq11}-\eqref{eq15} 10,000 times. When we apply the proposed algorithm in practice, we can expect faster convergence by introducing appropriate feasibility tolerances $\epsilon_\text{pri}$ and $\epsilon_\text{dual}$, and setting $\|r^{(k)}\| \leq \epsilon_\text{pri}$ and $\|s^{(k)}\| \leq \epsilon_\text{dual}$ as the termination conditions.
\begin{equation}
\label{eq30}
\begin{aligned}
  r^{(k+1)}
  & =
  M u^{(k+1)} - z^{(k+1)} - n \\
  s^{(k+1)}
  & =
  \rho M^T\left( z^{(k)} - z^{(k+1)} \right) 
\end{aligned}
\end{equation}

Figure \ref{fig9} and Figure \ref{fig10} show the miss distance $\|l_N\|$ and the impact angle errors $|\angle v_N - \theta_f|$ for different initial positions of the target; see the engagement scenario in Figure \ref{fig1}. The results are computed with fixed $\rho$ and fixed number of iterations for simplicity. which involves small guidance error. Note that the miss distance and the impact angle errors from the OGL solutions increase rapidly as the crosstrack of the target's initial position increases and the range decreases, whereas the errors are significantly smaller for the proposed approach.

\section{Concluding remarks}

We introduced a computational technique for optimal impact angle guidance problems, that directly handles nonconvex constraints arising from the kinematic constraints.
In general, the linearization techniques used for solving nonconvex problems are difficult to apply in large divert engagement geometries where linearization assumption is not valid and the computed solutions exhibits significant suboptimality. 
This paper solves these problems by using the first-order method with direct Euclidean projection onto the nonconvex sets. We present an efficient algorithm for computing orthogonal projection onto the nonconvex set described by an angular constraints between state vectors.

We verified the proposed approach through a series of numerical numerical simulations. The proposed technique was shown to obtain optimal trajectories that satisfy nonconvex constraints even in extreme maneuver cases, where the traditional impact angle guidance solution fails or returns significantly suboptimal solutions.

\bibliographystyle{IEEEtran}
\bibliography{refs}

\end{document}